\documentclass[10pt]{amsart}
\usepackage{latexsym, amsmath,amssymb}
\usepackage{amsmath}
\usepackage{amsfonts}
\usepackage{amssymb}
\usepackage{amscd}
\usepackage{amsthm}
\usepackage{amsbsy}
\usepackage{graphicx}
\usepackage{bm}

\renewcommand{\epsilon}{\varepsilon}
\newcommand{\newsection}[1]
{\subsection{#1}\setcounter{theorem}{0} \setcounter{equation}{0}
\par\noindent}

\newtheorem{theorem}{Theorem}

\newtheorem{lemma}[theorem]{Lemma}
\newtheorem{corr}[theorem]{Corollary}

\newtheorem{proposition}[theorem]{Proposition}
\newtheorem{deff}[theorem]{Definition}

\newcommand{\bth}{\begin{theorem}}
\newcommand{\ble}{\begin{lemma}}
\newcommand{\bcor}{\begin{corr}}

\newcommand{\bdeff}{\begin{deff}}

\newcommand{\bprop}{\begin{proposition}}
\newcommand{\ele}{\end{lemma}}
\newcommand{\ecor}{\end{corr}}
\newcommand{\edeff}{\end{deff}}

\newcommand{\eprop}{\end{proposition}}

\newcommand{\cd}{\, \cdot\, }

\newcommand{\la}{\lambda}

\newcommand{\supp}{\text{supp }}
\renewcommand{\Pi}{\varPi}

\renewcommand{\epsilon}{\varepsilon}

\newcommand{\Rt}{{\Bbb R}^2}

\newcommand{\parital}{\partial}
\newcommand{\tidle}{\tilde}

\newcommand{\Aut}{\text{Aut}(p)}

\begin{document}


\title[On restriction estimates and $L^4$-bounds]
{On eigenfunction restriction estimates and $L^4$-bounds for compact surfaces 
with nonpositive
curvature}
\thanks{The authors were supported in part by the NSF.  Some of this research was carried
out while the first author was visiting Zhejiang University in Hangzhou, China, and he wishes
to thank his colleagues there for their kindness. We also wish to thank W. Minicozzi III for many
very helpful discussions.}

\author{Christopher D. Sogge}
\author{Steve Zelditch}
\address{Johns Hopkins University, Baltimore, MD}
\address{Northwestern University, Evanston, IL}

\maketitle

\begin{abstract}
If $(M,g)$ be a two-dimensional compact
boundaryless Riemannian manifold with nonpostive curvature,
then we shall give improved estimates for the $L^2$-norms of
the restrictions of eigenfunctions to unit-length geodesics, compared to 
the general results of Burq, G\'erard and Tzvetkov
\cite{burq}.  By earlier results of 
Bourgain \cite{bourgainef} and the first author \cite{Sokakeya},
they are equivalent to improvements of the general $L^p$-estimates
in \cite{soggeest} for $n=2$ and $2<p<6$.  The proof uses 
the fact that the exponential map from any point in $x_0\in M$ is a universal
covering map from $\Rt \simeq T_{x_0}M$ to $M$ (the Cartan-Hadamard- von Mangolt theorem), which
allows us to lift the necessary calculations up to the universal cover
$(\Rt, \tilde g)$ where $\tilde g$ is the pullback of $g$ via the exponential
map.  We then prove the main estimates by using the Hadamard
parametrix for the wave equation on $(\Rt, \tilde g)$ and the fact that
the classical comparison theorem of G\"unther \cite{Gu} for the volume element 
in spaces of nonpositive curvature
gives us desirable bounds for the principal coefficient of the
Hadamard parametrix, allowing us to prove our main result.
\end{abstract}

\newsection{Introduction}

Let $(M,g)$ be a compact two-dimensional Riemannian manifold without boundary.  We shall assume
throughout that the curvature of $(M,g)$ is everywhere nonpositive.  If $\Delta_g$ is the
Laplace-Beltrami operator associated with the metric $g$, then we are concerned with
certain size estimates for the eigenfunctions
$$-\Delta_g e_\lambda(x)=\lambda^2 e_\lambda(x), \quad x\in M.$$
Thus we are normalizing things so that $e_\la$ is an eigenfunction of the first order
operator $\sqrt{-\Delta_g}$ with eigenvalue 
$\la$.
If $e_\la$ is also normalized to have $L^2$-norm one, we are interested
in various size estimates for the $e_\la$ which are related to how concentrated they
may be along geodesics.  If $\Pi$ denotes the space of all unit-length geodesics in $M$
then our main result is the following ``restriction theorem'' for this problem.

\begin{theorem}\label{theorem1}  Assume that $(M,g)$ is as above.  Then given
$\varepsilon>0$ there is a $\la(\varepsilon)<\infty$ so that
\begin{equation}\label{1.1}
\sup_{\gamma\in \Pi} \left(\, \int_\gamma |e_\la|^2 \, ds \, \right)^{1/2}
\le \varepsilon \la^{\frac14} \|e_\la\|_{L^2(M)}, \quad \la >\la(\varepsilon),
\end{equation}
with $ds$ denoting arc-length measure on $\gamma$, and $L^2(M)$ being
the Lebesgue space with respect to the volume element $dV_g$ for $(M,g)$.
\end{theorem}

Earlier, Burq, G\'erard and Tzvetkov \cite{burq} showed that for {\it any}
2-dimensional compact boundaryless Riemannian manifold one has
\begin{equation}\label{1.2}
\left(\, \int_\gamma |e_\la|^2 \, ds \, \right)^{1/2}
\le C \la^{\frac14} \|e_\la\|_{L^2(M)},
\end{equation}
with $C$ independent of $\gamma\in \Pi$.  
The first such estimates were somewhat weaker ones of Reznikov~\cite{rez} for
hyperbolic surfaces, which inspired this current line of research.
The estimate \eqref{1.2} is sharp for the round
sphere $S^2$ because of the highest weight spherical harmonics
(see \cite{burq}, \cite{Sokakeya}).  Burq, G\'erard and Tzvetkov \cite{burq}
also showed that
$$\left(\, \int_\gamma |e_\la|^4\, ds \, \right)^{1/4}\le C\la^{\frac14}\|e_\la\|_{L^2(M)},
\quad \gamma \in \Pi,$$
and so by interpolating with this result and \eqref{1.1} one concludes that
when $M$ has nonpositive curvature
$\sup_{\gamma\in \Pi}\|e_\la\|_{L^p(\gamma)}/\|e_\la\|_{L^2(M)}=o(\la^{\frac14})$
for $2\le p<4$.  An interesting but potentially difficult problem would be to show that
this remains true under this hypothesis for the endpoint $p=4$.

Theorem~\ref{theorem1} is related to certain $L^p$-estimates for eigenfunctions.  In
\cite{soggeest} the first author proved that for any compact Riemannian manifold of dimension
$2$ one has for $\la\ge1$,
\begin{equation}\label{1.3}
\|e_\la\|_{L^p(M)}\le C\la^{\frac12 (\frac12 -\frac1p)} \|e_\la\|_{L^2(M)},
\quad 2\le p\le 6,
\end{equation}
and
\begin{equation}\label{1.4}
\|e_\la\|_{L^p(M)}\le C\la^{2(\frac12-\frac1p)-\frac12}\|e_\la\|_{L^2(M)},
\quad  6\le p\le \infty.
\end{equation}
These estimates are also sharp for the round sphere $S^2$ (see \cite{sph}).  The first 
estimate, \eqref{1.3}, is sharp because of the highest weight spherical harmonics,
and thus, like \eqref{1.1} or \eqref{1.2}, it measures concentration of eigenfunction
mass along geodesics.  The second estimate, \eqref{1.4}, is sharp due to the
zonal functions on $S^2$, which concentrate at points.  The sharp variants of
\eqref{1.3} and \eqref{1.4} (with different exponents) for manifolds with boundary were obtained by H. Smith and
the first author in \cite{ss}, and it would be interesting to obtain analogues of 
the results in the present paper for this setting, but this appears to be difficult.

In the last decade there have been several results showing that, for typical $(M,g)$,
\eqref{1.4} can be improved for $p>6$ (see \cite{stz}, \cite{soggezelditch}) to bounds
of the form $\|e_\la\|_{L^p(M)}/\|e_\la\|_{L^2(M)}=o(\la^{\frac12 (\frac12 -\frac1p)})$
for fixed $p>6$.  Recently, Hassell and Tacey~\cite{HT}, following
B\'erard's \cite{Berard} earlier estimate for $p=\infty$, showed that
for fixed $p>6$ this ratio is $O(\la^{2(\frac12-\frac1p)-\frac12}/\sqrt{\log \la})$, which
influenced the present work. Also, in \cite{sz2} the authors showed that if the geodesic flow
is ergodic, which is automatically the case if the curvature of $M$ is negative,  then \eqref{1.1} holds for a
density one sequence of eigenfunctions.

Except for some special cases of an arithmetic nature (e.g. Zygmund~\cite{zygmund}
or Spinu~\cite{Spinu}) there have been few cases showing that \eqref{1.3}
can be improved for Lebesgue exponents with $2<p<6$.  In \cite{Sokakeya}, using in part results from Bourgain~\cite{bourgainef},
it was shown that
$$\|e_\la\|_{L^p(M)}/\|e_\la\|_{L^2(M)}=o(\la^{\frac12 (\frac12 -\frac1p)})$$ for some
$2<p<6$ if and only if 
$$\sup_{\gamma\in \Pi}\|e_\la\|_{L^2(\gamma)}/\|e_\la\|_{L^2(M)}
=o(\la^{\frac14}).$$  Thus, we have the following corollary to Theorem~\ref{theorem1}.

\begin{corr}\label{corollary2} As above, let $(M,g)$ be a compact $2$-dimensional
manifold with nonpositive curvature.  Then, if $\varepsilon>0$ and
$2<p<6$ are fixed there is a $\la(\varepsilon,p)<\infty$ so that
$$\|e_\la\|_{L^p(M)}\le \varepsilon
\la^{\frac12 (\frac12-\frac1p)}\|e_\la\|_{L^2(M)}, \quad \la >\la(\varepsilon,p).$$
\end{corr}

We remark that an interesting open problem would be to obtain this type of result for
the case of $p=6$.  It is valid for the standard torus ${\mathbb T}^2={\mathbb R}^2/{\mathbb Z}^2$ since Zygmund~\cite{zygmund} showed that there one has $\|e_\la\|_{L^4({\mathbb T}^2)}/\|e_\la\|_{L^2({\mathbb T}^2)}=O(1)$ and the classical theorem of Gauss about lattice points in the
plane yields $\|e_\la\|_{L^\infty({\mathbb T}^2)}/\|e_\la\|_{L^2({\mathbb T}^2)}=O(\la^{\frac14})$.
Since $p=6$ is the exponent for which concentration at points and concentration
along geodesics are both relevant, proving a general result along the lines
of Corollary~\ref{corollary2} would presumably have to take into account both of these
phenomena.  One expects, though, such a result for $p=6$ should be valid when $M$
has negative curvature.  This result seems to be intimately related to the problem of
trying to determine when one has the endpoint improvement for the
restriction problem, i.e., $\sup_{\gamma\in \Pi}\|e_\la\|_{L^4(\gamma)}/\|e_\la\|_{L^2(M)}
=o(\la^{\frac14})$.

In \cite{Sokakeya} the first author showed that if $\gamma_0\in \Pi$ is not part of a periodic geodesic
then $$\|e_\la\|_{L^2(\gamma_0)}/\|e_\la\|_{L^2(M)}=o(\la^{\frac14}).$$  The proof involved
an estimate involving the wave equation associated with $\Delta_g$ and a bit of
microlocal (wavefront) analysis.  The main step in proving Theorem~\ref{theorem1}
is to see that this remains valid as well if $\gamma_0$ is part of a periodic orbit under
the above curvature assumptions.  We shall be able to do this by lifting the wave equation
for $(M,g)$ up to the corresponding one for its universal cover, which
by a classical theorem of Hadamard~\cite{Had98} and von Mangolt~\cite{VM}, is $(\Rt, \tilde g)$,
with the metric $\tilde g$ being the pullback of $g$ via a covering map,
which can be taken to be $\exp_{x_0}$ for any $x_0\in M$.  By
identifying solutions of wave equations for $(M,g)$ with ``periodic'' ones for
$(\Rt, \tilde g)$ we are able to obtain the necessary bounds using a bit of wavefront
analysis and the Hadamard parametrix for $(\Rt, \tilde g)$.  Fortunately for us,
by a classical volume comparison theorem of G\"unther~\cite{Gu}, the leading
coefficient of the Hadamard parametrix has favorable size estimates under our
curvature assumptions.  (It is easy to see that the contribution of the lower order
terms in the Hadamard parametrix to \eqref{1.1} are straightforward to handle.)

\newsection{Proof of geodesic restriction bounds}

Since the space of all unit-length geodesics is compact, in order to
prove \eqref{1.1}, it suffices to show that, given $\gamma_0\in \Pi$ and
$\varepsilon>0$, one can find a neighborhood ${\mathcal N}(\gamma_0,\varepsilon)$
of $\gamma_0$ in $\Pi$ and a number $\la(\gamma_0,\varepsilon)$ so that
\begin{equation}\label{2.1}
\int_{\gamma}|e_\la|^2 \, ds \le \varepsilon \la^{\frac12}\|e_\la\|^2_{L^2(M)},
\quad \gamma\in {\mathcal N}(\gamma_0,\varepsilon), \, \, \, \la>\la(\gamma_0,\varepsilon).
\end{equation}
In proving this we may assume that the injectivity radius of $(M,g)$ is ten or more.  We
recall also that, given $x_0\in M$, the exponential map at $x_0$,
$\exp_{x_0}: T_{x_0}M\simeq \Rt \to M$ is a universal covering map.  We shall take 
$x_0$ to be the midpoint of our unit-length geodesic $\gamma_0$.  We also shall
work in geodesic polar coordinates about $x_0$.

If $\tilde g$ is the pullback to $\Rt$  of the metric $g$ via the covering map
then $(\Rt,\tilde g)$ is a Riemannian universal cover of $(M,g)$.  Like
$(M,g)$ it also has
nonpositive curvature.  Additionally, rays $t\to t(\cos \theta, \sin \theta)$, $t\ge0$, through the origin are geodesics for $\tilde g$.  Such a ray is the lift of the unit speed geodesic starting
at $x_0$, which in our local coordinate system has the initial tangent vector
$(\cos \theta, \sin \theta)$.  Note that in these coordinates vanishing at $x_0$, 
$t\to t(\cos\theta, \sin\theta)$, $|t|\le 10$ are also geodesics for $g$.  We may assume further
that we have
\begin{equation}\label{2.2}
\gamma_0 =\{(t,0): \, -\frac12 \le t\le \frac12 \}.
\end{equation}

To prove \eqref{2.1} it will be convenient to fix a real-valued  even function $\chi\in {\mathcal S}
({\mathcal R})$ having the property that $\chi(0)=1$ and $\Hat \chi(t)=0$, $|t|\ge \frac14$,
where $\Hat \chi$ denotes the Fourier transform of $\chi$.  We then have that for
$T>0$
$$\chi(T(\sqrt{-\Delta_g}-\la))e_\la = e_\la,$$
and, therefore, to prove \eqref{2.1}, it suffices to show that if $T$ is large and
fixed then there is a neighborhood
\footnote{We can use the topology of $S^*M$ to define these neighborhoods, since every $\gamma\in \Pi$ can be uniquely identified with an element $(y,\xi)\in S^*M$, with $y$ being
the midpoint of $\gamma$ and $\xi$ being the direction of $\gamma$ at $y$.}
 ${\mathcal N}={\mathcal N}(\gamma_0,T)$
of $\gamma_0$ so that
\begin{equation}\label{2.3}
\int_\gamma \bigl| \,
\chi(T(\sqrt{-\Delta_g}-\la))f\, \bigr|^2 \, ds
\le CT^{-1}\la^{\frac12}\|f\|_{L^2(M)}^2 +
C'_{T,{\mathcal N}}\|f\|_{L^2(M)}^2, \quad \gamma \in {\mathcal N},
\end{equation}
where $C$ (but not $C'_{T,{\mathcal N}}$) is a uniform constant
depending on $(M,g)$ but independent of $T$ and ${\mathcal N}$.

To prove \eqref{2.3}, we shall be able to use the wave equation as
\begin{multline}\label{2.4}
\chi(T(\sqrt{-\Delta_g}-\la))f=\frac1{2\pi T}
\int_{{\mathbb R}}\Hat \chi(t/T) e^{-it\la} e^{it\sqrt{-\Delta_g}}f\, dt
\\
=\frac1{\pi T}\int_{-T/4}^{T/4}\Hat \chi(t/T) e^{-it\la}\cos t\sqrt{-\Delta_g}f \, dt
+\chi(T(\sqrt{-\Delta_g}+\la))f,
\end{multline}
using the fact that $\Hat \chi(t)$ is even and supported in $|t|\le \frac14$.  Since
the kernel of the last term satisfies
\begin{equation}\label{2.5}
|\partial^\alpha_{x,y}\chi(T(\sqrt{-\Delta_g}+\la))(x,y)|\le C_{T,N}\la^{-N}
\end{equation}
for any $N$ in compact subsets of any local coordinate system, to prove
\eqref{2.3} it suffices to show that
\begin{equation}\label{2.6}
\int_{\gamma}
\left| \, \frac1{\pi T}\int_{-T/4}^{T/4} \Hat \chi(t/T) e^{-it\la} \cos t\sqrt{-\Delta_g}f\, dt \, 
\right|^2 \, ds 
\le \bigl(\, CT^{-1}\la^{\frac12} +C'_{T,\mathcal{N}}\, \bigr) \|f\|_{L^2(M)}, \quad
\gamma\in {\mathcal N}(\gamma_0,T).
\end{equation}

If $\gamma_0$ is not part of a periodic geodesic of period $\le T$, then we can easily prove 
\eqref{2.6} just by using wavefront analysis and arguments that are similar to the
proof of the Duistermaat-Guillemin theorem \cite{DG}.  This was done in
\cite{Sokakeya}, but we shall repeat the argument here for the sake of
completeness and since it motivates what is needed to handle the argument
when $\gamma_0$ is a portion of  a periodic geodesic of period $\le T$.

To handle the latter case we shall exploit the relationship between solutions
of the wave equation on $(M,g)$ of the form
\begin{equation}\label{2.7}
\begin{cases}
(\partial^2_t-\Delta_g)u(t,x)=0, \quad (t,x)\in {\mathbb R}_+\times M
\\
u(0,\cd)=f, \, \, \partial_t u(0,\cd)=0,
\end{cases}
\end{equation}
and certain ones on $({\mathbb R},\tilde g)$
\begin{equation}\label{2.8}
\begin{cases}
(\partial^2_t-\Delta_{\tilde g})\tilde u(t,\tilde x), \quad (t,\tilde x)\in {\mathbb R}_+\times
\Rt
\\
\tilde u(0,\cd)=\tilde f, \, \, \partial_t \tilde u(0,\cd)=0.
\end{cases}
\end{equation}
Note that
$u(t,x)=\bigl(\cos (t\sqrt{-\Delta_g})f\bigr)(x)$
is the solution of \eqref{2.7}.

To describe the relationship between the two equations
we shall use the deck transformations
associated with our universal covering map
\begin{equation}\label{2.9}
p=\exp_{x_0}: \, \Rt \to M.
\end{equation}
Recall that an automorphism for $(\Rt, \tilde g)$,  $\alpha: \, \Rt \to \Rt$ ,is a deck transformation if
$$p\circ \alpha = p.$$
In this case we shall write $\alpha \in \Aut$.  In the case
where ${\mathbb T}^2$ is the standard two-torus, each
$\alpha$ would just be translation in $\Rt$ with respect
to some $j\in {\mathbb Z}^2$.  Motivated by this if 
$\tilde x\in\Rt$ and $\alpha\in \Aut$, let us call 
$\alpha(\tilde x)$ the translate of $\tidle x$ by $\alpha$.
then we recall a set $D\subset \Rt$ is called a fundamental domain
of our universal covering $p$ if every point in $\Rt$ is the
translate of exactly one point in $D$.  Of course there are
infinitely many fundamental domains, but we may assume that
ours is relatively compact, connected and contains the ball of radius
2 centered at the origin in view of our assumption about
the injectivity radius of $(M,g)$.  We can then think of our unit geodesic
$\gamma_0=\{(t,0): \, |t|\le \frac12\}$ (written in geodesic polar coordinates as
above) both as one in $(M,g)$ and one in the fundamental domain which is of 
the same form.  Likewise, a function $f(x)$ on $M$ is uniquely identified
by one $f_D(\tilde x)$ on $D$ if we set $f_D(\tilde x)=f(x)$, where
$\tilde x$ is the unique point in $D\cap p^{-1}(x)$.  Using $f_D$ we can
define a ``periodic extension'', $\tilde f$, of $f$ to $\Rt$ by defining 
$\tilde f (\tilde y)$ to be equal to $f_D(\tilde x)$ if $\tilde x=\tilde y$ modulo
$\Aut$, i.e. if $(\tilde x,\alpha)\in D\times \Aut$ are the unique pair so that
$\tilde y=\alpha(\tilde x)$.  Note then that $\tilde f$ is periodic with respect
to $\Aut$ since we necessarily have that
$\tilde f(\tilde x)=\tilde f(\alpha(\tilde x))$ for every $\alpha \in \Aut$.

We can now describe the relationship between the wave equations
\eqref{2.7} and \eqref{2.8}.  First, if $(f(x),0)$ is the Cauchy
data in \eqref{2.7} and $(\tilde f(\tilde x),0)$ is the
periodic extension to $(\Rt, \tilde g)$, then the solution $\tilde u(t,\tilde x)$
to \eqref{2.8} must also be a periodic function of $\tilde x$ since 
$\tilde g$ is the pullback of $g$ via $p$ and $p=p\circ \alpha$.  As a 
result, we have that the solution to \eqref{2.7} must satisfy
$u(t,x)=\tilde u(t,\tilde x)$ if $\tilde x\in D$ and $p(\tilde x)=x$.  Another
way of saying this is that if $\tilde f$ is the pullback of $f$ via $p$ and
$t$ is fixed then $\tilde u(t,\cd)$ solving \eqref{2.8} must be the pullback
of $u(t,\cd)$ in \eqref{2.7}.  Thus, periodic solutions to \eqref{2.8}
correspond uniquely to solutions of \eqref{2.7}.  In other words,
we have the important formula for the wave kernels
\begin{equation}\label{2.10}
\bigl(\cos(t\sqrt{-\Delta_g})(x,y)=\sum_{\alpha\in \Aut}
\bigl(\cos (t\sqrt{-\Delta_{\tilde g}}\bigr)(\tidle x, \alpha (\tilde y)), 
\end{equation}
if $\tilde x$ and $\tilde y$ are the unique points in $D$ for which
$p(\tilde x)=x\quad \text{and } \, \, p(\tilde y)=y$.

Note that the sum in \eqref{2.10} only has finitely many nonzero terms
for a given $(x,y,t)$ since, by the finite propagation speed for 
$\square_{\tilde g}=\parital_t^2-\Delta_{\tilde g}$, the summands in the
in the right all vanish when $d_{\tilde g}(\tilde x,\alpha(\tilde y))>t$.  
For instance, if $x=y=x_0$ the number of nontrivial terms would equal
the cardinality of $p^{-1}(x_0)\cap \{\tilde x\in \Rt: \, |\tilde x|\le t\}$
where $|\tilde x|$ denotes the Euclidean length, due to the fact
that $d_{\tilde g}(0,\tilde x)=|\tilde x|$.  Despite this, the number
of nontrivial terms will grow exponentially in $t$ if the curvature
is bounded from above by a fixed negative constant.

To see this, let us review one last thing before focusing more closely
on the proof of our restriction-estimate.  As we shall see, even though
there can be an exponentially growing number of nontrivial terms in the right hand side of
\eqref{2.10}, which could create havoc for our proofs if we are not
careful, this turns out to be related to something that will actually
be beneficial for our calculations.

These facts are related to the fact that
in  the geodesic polar coordinates we are
using, $(t \cos \theta, t \sin \theta)$, $t>0$, $\theta\in (-\pi,\pi]$,
for $(\Rt, \tilde g)$,
the metric $\tilde g$ takes the form
\begin{equation}\label{2.11}
ds^2=dt^2+\mathcal{A}^2(t,\xi)\, d\theta^2,
\end{equation}
where we may assume that $\mathcal{A}(t,\theta)>0$ for $t>0$.
Consequently, the volume element in these coordinates is
given by
\begin{equation}\label{2.12}
dV_g(t,\theta)=\mathcal{A}(t,\theta)\, dt d\theta,\end{equation}
and by G\"unther's \cite{Gu} comparison theorem if the curvature
of $(M,g)$ and hence that of $(\Rt, \tilde g)$ is nonpositive, we have
\begin{equation}\label{2.13}\mathcal{A}(t,\theta)\ge t.
\end{equation}
Furthermore, if one assumes that the curvature is $\le - \kappa^2$,
with $\kappa>0$ then one has 
\begin{equation}\label{2.14}
\mathcal{A}(t,\theta)\ge \frac1\kappa \sinh (\kappa t).\end{equation}
Since the volume element for two-dimensional Euclidean space in polar coordinates
is $t \, dt d\theta$ and that of the hyperbolic plane with constant curvature
$-\kappa^2$ is $\frac1\kappa \sinh(\kappa t) \, dt d\theta$, G\"unther's 
volume comparison theorem
says that in geodesic polar coordinates the volume element for spaces
of nonpositive curvature is at least that of $\Rt$ with the flat metric, while
if the curvature is bounded above by $-\kappa^2$ the volume element is
at least that of the hyperbolic plane of constant curvature $-\kappa^2$.
In the latter case, as we warned, the number of nontrivial terms
in the sum in the right side of \eqref{2.10} will be at least bounded
below by a multiple of $e^{\kappa t}$ as $t\to +\infty$.

Let us now turn to the proof of \eqref{2.6} and hence Theorem~\ref{theorem1}.  
Given $\gamma\in\Pi$ we let  $T^*\gamma\subset
T^*M$ and $S^*\gamma\subset S^*M$ be the cotangent and unit cotangent
bundles over $\gamma$, respectively.  Thus, if $(x,\xi)\in T^*\gamma$ then
$\xi_\sharp$ is a tangent vector  to $\gamma$ at $x$ if $T^*M\ni\xi\to \xi_\sharp\in TM$ is the
standard musical isomorphism, which, in local coordinates, sends $\xi=(\xi_1,\xi_2)
\in T^*_xM$ to $\xi_\sharp=(\xi_\sharp^1,\xi_\sharp^2)$ with
$\xi_\sharp^j=\sum_k g^{jk}(x)\xi_k$.  Then if $\Phi_t: S^*M\to S^*M$ denotes geodesic flow in the unit cotangent bundle
over $M$, and $(x,\xi)\in S^*\gamma$
we let
$L(x,\xi)$ be the minimal $t>0$ so that $\Phi_t(x,\xi)=(x,\xi)$ and define it to be $+\infty$
if no such time $t$ exists.  Then if $\gamma$ is not part of a periodic geodesic this
quantity is $+\infty$ on $S^*\gamma$, and if it is then it is constant on $S^*\gamma$ and equal to the
minimal period of the geodesic, $\ell(\gamma)$ 
(which must be larger than ten because of our
assumptions).  Note also that $L(x,\xi)$ can also be thought of
as a function on $S^*M$, and that, in this case, it is lower semicontinuous.

Recall
that we are working in geodesic polar coordinates vanishing at $x_0$, the midpoint
of $\gamma_0$, and that $\gamma_0$ is of the form \eqref{2.2} in these
coordinates.  Let us choose $\beta\in C^\infty_0({\mathbb R})$ equal to one
on $[-\frac34, \frac34]$ but $0$ outside $[-1,1]$.  We then let $b_\varepsilon(x,D)$
and $B_\varepsilon(x,D)$ be zero-order pseudodifferential operators which in
the above local coordinates have symbols
$$b_\varepsilon(x,\xi)=\beta(|x|)\beta(\xi_2/\varepsilon |\xi|), \quad \text{and } \, B_\varepsilon(x,\xi)=
\beta(|x|)(1-\beta(\xi_2/\varepsilon|\xi|)),$$
respectively.

Our first claim is that if $\varepsilon>0$ and $\gamma\in \Pi$ are fixed, then we can find
a neighborhood ${\mathcal N}(\gamma_0,\varepsilon)$ of $\gamma_0$ so that
\begin{equation}\label{2.15}
\int_{-T/4}^{T/4}\int_\gamma\left| \, B_\varepsilon \circ \cos (t\sqrt{-\Delta_g})f\, \right|^2 \, ds dt
\le C_{T,\varepsilon}\|f\|_{L^2(M)}^2, \quad \gamma \in \mathcal{N}(\gamma_0,\varepsilon),
\end{equation}
which, by an application of the Schwartz inequality, would yield part of \eqref{2.6}, namely,
\begin{equation}\label{2.16}
\int_\gamma \Bigl| \, \frac1{\pi T}\int_{-T/4}^{T/4} \Hat \chi(t/T) e^{-i\la t}B_\varepsilon
\circ \cos(t\sqrt{-\Delta_g})f\, dt \, \Bigr|^2 \, ds
\le C'_{T,\varepsilon}\|f\|_{L^2(M)}^2, \quad \gamma \in {\mathcal N}(\gamma_0,\varepsilon).
\end{equation}
If $R_\gamma$ denotes the restriction to $\gamma\in \Pi$, then \eqref{2.15} follows from the
fact that the operator
$$f\to R_\gamma (A\circ \cos(t\sqrt{-\Delta_g})f),$$
regarded as an operator from $C^\infty(M)\to C^\infty(\gamma \times [-T/4,T/4])$, is a Fourier
integral operator of order zero which is locally a canonical graph (i.e., nondegenerate)
if $\text{supp }A(x,\xi)\cap S^*\gamma=\emptyset$, and hence a bounded operator
from $L^2(M)$ to $L^2(\gamma \times [-T/4,T/4])$.  since $B_\varepsilon(x,\xi)$ vanishes
on a neighborhood of $S^*\gamma_0$, we conclude that this is the case $A=B_\varepsilon$ for $\gamma\in
\Pi$ close to $\gamma_0$, which gives us \eqref{2.15}.  The $L^2$-boundedness of 
nondegenerate Fourier integrals is a theorem of H\"ormander~\cite{hormander2}, while
the observation about $R_\gamma (A\circ \cos(t\sqrt{-\Delta_g}))$ is one of  
Tataru~\cite{tataru}.  It is also easy to check the latter, because, for fixed $t$,
$e^{it\sqrt{-\Delta_g}}: C^\infty(M)\to C^\infty(M)$ is a nondegenerate Fourier integral
operator, and, therefore, one needs only to verify the assertion when $t=0$, in which
case it is an easy calculation using any parametrix for the half-wave operator.

The estimate \eqref{2.16} holds for any $\gamma_0\in \Pi$.  Let us now argue
that if $\ell(\gamma_0)$, the period of $\gamma_0$, is larger than $T$
or if $\gamma_0$ is not part of a periodic geodesic, then we have
also have favorable bounds if $B_\varepsilon$ is replaced by $b_\varepsilon$,
with $\varepsilon>0$ sufficiently small.  To do this, we recall that the wave front
set of the kernel of
$b_\varepsilon \circ \cos(t\sqrt{-\Delta_g}) \circ b_\varepsilon^*$ is contained in
\begin{equation}\label{2.17}
\bigl\{(x,t,\xi,\tau; y, -\eta): \, \Phi_{\pm t}(x,\xi)=(y,\eta), \, \, \tau^2 =\sum g^{jk}(x)\xi_j\xi_k, 
\, \, 
(x,\xi), (y,\eta)\in 
\supp b_\varepsilon\bigr\}.
\end{equation}
To exploit this, let $W_\gamma$ be the operator
\begin{equation}\label{2.17'}
W_\gamma f=R_\gamma\Bigl(\, \frac1{\pi T}\int_{-T/4}^{T/4}
\Hat \chi(t/T) e^{-i\la t} b_\varepsilon \circ \cos (t\sqrt{-\Delta_g})f \, dt \, \Bigr).\end{equation}
Our goal then is to show, that under the present assumption that
$\ell(\gamma_0)> T$,
$$\|W_\gamma\|_{L^2(M)\to L^2(\gamma)} \le CT^{-\frac12}\la^{\frac14}+C_{T,b_\varepsilon}$$
for $\gamma\in \Pi$ belonging to some neighborhood ${\mathcal N}(\gamma_0,T,\varepsilon)$
of $\gamma_0$.  This is equivalent to showing that the dual operator $W_\gamma^*: L^2(\gamma)
\to L^2(M)$ with the same norm, and since
$$\|W_\gamma^* g\|_{L^2(M)}^2 =
\int_\gamma W_\gamma W_\gamma^*g \, \overline{g}\, ds \le \|W_\gamma W_\gamma^*g\|_{L^2(\gamma)}\|g\|_{L^2(\gamma)},$$
we would be done if we could show that
\begin{equation}\label{2.18}
\|W_\gamma W_\gamma^*g\|_{L^2(\gamma)}\le \Bigl( \, CT^{-1}\la^{\frac12}+C_{T,b_\varepsilon}\, \Bigr)
\|g\|_{L^2(\gamma)}.
\end{equation}
But, by Euler's formula, the kernel of $4W_\gamma W_\gamma^*$ is $K|_{\gamma \times \gamma}$, where
$K(x,y)$, $x,y\in M$ is the kernel of the operator
$$b_\varepsilon \circ \rho(T(\sqrt{-\Delta_g}-\la))\circ b_\varepsilon^*
+b_\varepsilon \circ \rho(T(\sqrt{-\Delta_g}+\la))\circ b_\varepsilon^*
+2b_\varepsilon \circ \chi(T(\sqrt{-\Delta_g}-\la)) \chi(T(\sqrt{-\Delta_g}+\la))\circ b_\varepsilon^*
,$$
if $\rho(\tau)=(\chi(\tau))^2$.
The last two terms satisfy bounds like those in \eqref{2.5} (with constant depending
on $T$ and $b_\varepsilon$), and the first term is
\begin{equation}\label{2.19}
\frac1{\pi T} \int_{-T/2}^{T/2}\Hat \rho (t/T) e^{-i\la t}
\bigl( \, b_\varepsilon \circ \cos(t\sqrt{-\Delta_g})\circ b^*_\varepsilon\, \bigr)(x,y) \, dt.
\end{equation}
We are using the fact that $\Hat \rho =  \Hat \chi * \Hat \chi$ is 
supported in $[-\frac12, \frac12]$.  In view of \eqref{2.17}, if $\varepsilon>0$ is sufficiently
small, since  we are assuming that $\ell(\gamma_0)> T$, it follows that we can find a neighborhood ${\mathcal N}$
of $\gamma_0$ in $M$ so that $(b_\varepsilon \circ \cos(t{\sqrt-\Delta_g})\circ b_\varepsilon^*)(x,y)$
is smooth on ${\mathcal N}\times {\mathcal N}$ when $t\ge 2$.  Thus, on ${\mathcal N}\times
{\mathcal N}$ the difference between \eqref{2.18} and
$$K(x,y)=
\frac1{\pi T} \int_{-T/2}^{T/2} \beta(t/5) \, 
 \rho(t/T) e^{-it\la} 
\Bigl(b_\varepsilon \circ \cos(t\sqrt{-\Delta_g})\circ b_\varepsilon^*\Bigr)(x,y) \, dt$$
is $O_{T,b_\varepsilon}(1)$.  But, by using the Hadamard parametrix
(see
below) one finds that
\begin{equation}\label{2.20}|K(x,y)|\le CT^{-1}\la^{\frac12}(d_g(x,y))^{-\frac12}+
C_{b_\varepsilon,T}
\,
(1+ \la \,  (1+\la d_g(x,y))^{-\frac32}), \quad
x,y\in {\mathcal N},\end{equation}
for some uniform constant $C$, which is independent of $\varepsilon$, $T$ and $\la$. 
Since, by Young's inequality, the integral operator with kernel $K|_{\gamma\times \gamma}$ is bounded from
$L^2(\gamma)\to L^2(\gamma)$ with norm bounded by $CT^{-1}\la^{\frac12}
+C_{b_\varepsilon, T}$ if $\gamma\subset {\mathcal N}$, we get \eqref{2.18},
which finishes the proof that \eqref{2.3} holds provided that $\ell(\gamma_0) >T$.

The above argument used the fact that if $\ell(\gamma_0)>T$, with
$T$ fixed, then if $\varepsilon>0$
is small enough and $(x,\xi)\in \supp b_\varepsilon$ with $x\in \gamma_0$ then
$\Phi_t(x,\xi)\notin \supp b_\varepsilon$ for $2<|t|\le T/2$. 
In effect, this allowed us to cut the effect of loops though $\gamma_0$ of its extension
of length $T$ from our main calculuation, since they were all transverse.
 If  $\gamma_0
\in \Pi$ is part of a periodic geodesic of period $\le T$, i.e., $\ell(\gamma_0)\le T$, then
this need not be true.  On the other hand, if $T$ is fixed and $(x,\xi)$ is as above,
then for sufficiently small $\varepsilon$ we will have
\begin{equation}
\label{2.21}
\Phi_{\pm t}(x,\xi)\notin \supp b_\varepsilon, \, \, \text{if } x\in \gamma_0, \, \,
\text{and \, } t\notin \bigcup_{j\in {\mathbb Z}} \bigl[ \, j \ell(\gamma_0) -2, j\ell(\gamma_0)+2 \, \bigr].
\end{equation}
Note that our assumption that the injectivity radius of $(M,g)$ is 10 or more implies
that 
$$\ell(\gamma_0)\ge 10.$$  

To exploit this, we shall use \eqref{2.10} which relates the wave kernel for $(M,g)$ with
the one for its universal cover using the covering map given by $p=\exp_{x_0}$
with $x_0$ being the midpoint of $\gamma_0$.  Note that the points $\alpha(0)$,
$\alpha\in \Aut$ exactly correspond to geodesic loops through $x_0$, with looping
time being equal to the distance from $\alpha(0)$ to the origin in ${\mathbb R}^2$.
Just a few of these correspond to smooth loops through $x_0$ along the periodic
geodesic containing $\gamma_0$.  Since we are assuming that we are working with
local coordinates on $(M,g)$ and global geodesic polar ones on $(\Rt,\tilde g)$ so that
$\gamma_0$ is of the form \eqref{2.2}, the automorphisms with this property are
exactly the $\alpha_j\in \Aut$, $j\in {\mathbb Z}$ for which
\begin{equation}\label{2.22}\alpha_j(0)= \bigl( \, j\ell(\gamma_0), \, 0\, \bigr).
\end{equation}
Note that $G_{\gamma_0}=\{\alpha_j\}_{j\in {\mathbb Z}}$ is a cyclic subgroup of $\Aut$ with generator
$\alpha_1$, which 
 is the stabilizer group for the lift of periodic geodesic containing $\gamma_0$.
Consequently,  we can choose $\varepsilon>0$ small enough 
and a neighborhood ${\mathcal N}$ of $\gamma_0$ in $M$
so that\footnote{We should point out that we are abusing notation a bit in \eqref{2.23}.  The last
factor denotes the kernel of the integral operator on $M$ 
with kernel $H(x,y)=(\cos t\sqrt{-\Delta_{\tilde g}})(\tilde x,\alpha(\tilde y))$
is composed on the left and right by $b_\varepsilon$ and $b_\varepsilon^*$, respectively, and as 
before we are identifying points $x$ in $M$ with their cousin $\tilde x$ in the fundamental
domain.  In the coordinate systems we are using, though, both are the same when we
are close to $\gamma_0$.}
\begin{equation}\label{2.23}
\bigl(b_\varepsilon \circ \cos(t\sqrt{-\Delta_{\tilde g}})\circ b_\varepsilon^*\bigr)(\tilde x, \alpha(\tilde y))
\in C^\infty\bigl({\mathcal N}\times {\mathcal N} \times 
\bigl[ \, j \ell(\gamma_0) -2, j\ell(\gamma_0)+2 \, \bigr]\bigr), 
\quad \text{if } \, \, \Aut\ni \alpha \notin G_{\gamma_0}.
\end{equation}
Therefore, by \eqref{2.21}--\eqref{2.23}, if we repeat the arguments that were used
to prove \eqref{2.18}, we conclude that we would have
\begin{equation}\label{2.24}
\int_\gamma \Bigl| \, \frac1{\pi T}\int_{-\infty}^{\infty} \Hat \chi(t/T) e^{-i\la t}b_\varepsilon
\circ \cos(t\sqrt{-\Delta_g})f\, dt \, \Bigr|^2 \, ds
\le \bigl(\, CT^{-\frac14}\la^{\frac14} + C_{T,b_\varepsilon}\, \bigr)^2\|f\|_{L^2(M)}^2, 
\gamma \in {\mathcal N}(\gamma_0,T),
\end{equation}
for some neighborhood ${\mathcal N}(\gamma_0,T)$ in $\Pi$, if we could show that
if the $\alpha_j$ are as in \eqref{2.22} and
\begin{multline}\label{2.25}
K(x,y)=
\frac1{\pi T} \sum_{\{j\in {\mathbb Z}_+: \, j t(\gamma_0)\le T/2\}}
\int_{-\infty}^{\infty}
\beta\bigl((s-j\ell(\gamma_0))/5\bigr) \Hat \rho(s/T) e^{-is\la}
\\
\times
\bigl( \, b_\varepsilon \circ \cos s\sqrt{-\Delta_{\tilde g}}\circ b^*_\varepsilon  \, \bigr)(\tidle x, 
\alpha_j(\tilde y))\, ds,
\end{multline}
then 
\begin{equation}\label{2.26}
|K(x,y)|
\le \bigl(\, CT^{-1}\la^{\frac12}(d_g(x,y))^{-\frac12}+T^{-\frac12}\la^{\frac12}+C_{T,b_\varepsilon}\, 
(1+ \la \, (1+\la d_g(x,y))^{-\frac32}) \, 
\bigr), \quad x,y\in {\mathcal N},
\end{equation}
with ${\mathcal N}$ being some neighborhood in $M$ of $\gamma_0$ (depending
on $T$).  The second term in the right side of this inequality did not occur in the previous
steps.  It comes from the terms in \eqref{2.25} with $j\ne 0$.  Also, the fact that \eqref{2.26}
yields \eqref{2.24} just follows from an application of Young's inequality.

To prove \eqref{2.26}, it suffices to see that we can find ${\mathcal N}$
as above so that 
\begin{multline}\label{2.27}
\Bigl| \int \beta\bigl((s-j\ell(\gamma_0))/5\bigr) \Hat \rho(s/T) e^{-is\la}
\bigl( \, b_\varepsilon \circ \cos s\sqrt{-\Delta_{\tilde g}}\circ b^*_\varepsilon  \, \bigr)(\tidle x, 
\alpha_j(\tilde y))\, ds\Bigr|
\\
\le 
C\la^{\frac12} 
  \, \bigl(\, \max\bigl\{d_g(\tilde x, \alpha_j(\tilde y)), 
e^{\kappa d_g(\tilde x, \alpha_j(\tilde y))}
\bigr\}\, \bigr)^{-\frac12}\,
+C_{T,b_\varepsilon}
, \, \, \, 
x,y\in {\mathcal N}, \,  \, \,  0\ne |j| \ell(\gamma_0)\le T,
\end{multline}
assuming that the curvature of $(M,g)$ is everywhere $\le-\kappa^2$, $\kappa\ge0$, while
for $j=0$, we have
\begin{multline}\label{2.28}
\int \beta(s/5) \Hat \rho(s/T) e^{-is\la}
\bigl( \, b_\varepsilon \circ \cos s\sqrt{-\Delta_{ g}}\circ b^*_\varepsilon  \, \bigr)(x, 
 y))\, ds\Bigr|
\\
\le C \la^{\frac12}(d_g(x,y))^{-\frac12} +C_{T,b_\varepsilon} \bigl(\, 1+\la\, (1+\la d_g(x,y))^{-\frac32}\, \bigr), \quad x,y\in {\mathcal N}.
\end{multline}
Note that $d_{\tilde g}(\tilde x, \alpha_j(\tilde y))\in [j \ell(\gamma_0)-1, 
j\ell(\gamma_0)+1]$ when $x, y\in 
\gamma_0$ and hence $d_{\tilde g}(\tilde x, \alpha_j(\tilde y))\ge |j|$ when 
 $x,y\in {\mathcal N}$ with ${\mathcal N}$ being a small
neighborhood of $\gamma_0$ in $M$.  We shall assume that this is the case in
what follows.
We then get \eqref{2.26} by summing over $j$.
(Observe that if the curvature is assumed to
be bounded below by a negative constant, we get something a bit stronger than
\eqref{2.26} where in the second term we may replace $T^{-\frac12}$ by $T^{-1}$.)

Both \eqref{2.27} and \eqref{2.28} are routine consequences of stationary phase
and the Hadamard parametrix for the wave equation.

To prove \eqref{2.28} let $\phi(x,y)$ denote geodesic normal coordinates of
$y$ about $x$.  Then if $|t|\le 5$, by the Hadamard parametrix (see \cite{Hor3} or
\cite{SHZ}) and the composition calculus for Fourier integral operators
(see Chapter 6 in \cite{soggebook})
\begin{equation}\label{2.29}
\bigl(b_\varepsilon\circ \cos (t\sqrt{-\Delta_g})\circ b^*_\varepsilon
\bigr)(x,y) = \sum_{\pm} \int_{\Rt} e^{i\phi(x,y)\cdot \xi \pm it|\xi|}a_\varepsilon(x,y,\xi)\, d\xi
+O_\varepsilon(1),
\end{equation}
where $a_\varepsilon \in S^0_{1,0}$ depends on $-\Delta_g$ and $b_\varepsilon$ but
satisfies
\begin{equation}\label{2.30}
|a_\varepsilon|\le C, \quad \text{and } \, \,
|\partial^\alpha_{x,y}\partial^\sigma_\xi a_\varepsilon|\le C_{\varepsilon \alpha\sigma}
(1+|\xi|)^{-|\sigma|}.
\end{equation}
The first constant is independent of $C$ and only depends on the size of the symbol
of $b_\varepsilon$, which is $\le \|\beta\|_{L^\infty({\mathbb R})}^4$.
Recall (see \cite{soggebook}) the following fact about the Fourier transform of 
a density times Lebesgue measure on the circle
$S^1=\{\Theta =(\cos \theta, \sin \theta)\}$,
\begin{equation}\label{2.31}
\int_0^{2\pi} e^{iw\cdot \Theta} a_\varepsilon (x,y,\Theta) \, d\theta
=|2\pi w|^{-\frac12} \sum_{\pm}e^{\pm i |w|}a_\varepsilon(x,y, \pm w) 
+O_\varepsilon(|w|^{-\frac32}), \quad |w|\ge 1,
\end{equation}
where the constants for the last term depend on the size of finitely many constants
in \eqref{2.31}.  Since $|\phi(x,y)|=d_g(x,y)$, if we combine \eqref{2.29} and
\eqref{2.30}, we find that, modulo a $O_\varepsilon(1)$ term, if
$\psi(s)=\beta(s/5)\Hat \rho(s/T)$, then when $d_g(x,y)\ge \la^{-1}$, the
quantity in \eqref{2.28} is the sum over $\pm$ of a fixed multiple of
\begin{multline*}
(d_g(x,y))^{-\frac12}
\int_0^\infty \bigl(\, \Hat \psi(\la-r)+\Hat \psi(\la+r)\, \bigr) e^{\pm i rd_g(x,y)}
a_\varepsilon(x,y, \pm r\phi(x,y)) r^{\frac12} \, dr
\\
+O_\varepsilon\bigl(\,(d_g(x,y))^{-\frac32} \int_0^\infty
( |\Hat \psi (\la-r)|+|\Hat \psi(\la+r)|) \, (1+r)^{-\frac12} \, dr\, \bigr).
\end{multline*}
By \eqref{2.30}, the first term in is $O(\|a_\varepsilon\|_\infty (\la d_g(x,y))^{-\frac12})$, 
since $|\Hat \psi(\tau)|\le C_N(1+|\tau|)^{-N}$ for any $N$.  Since the last term
is $O_\varepsilon(\la^{-1/2}(d_g(x,y))^{-\frac32})$,   we have established \eqref{2.28}
when $d_g(x,y)\ge \la^{-1}$.  The fact that it is also $O(\la)+O_\varepsilon(1)$ is a simple
consequence of \eqref{2.29} and \eqref{2.30}
which gives the bounds for $d_g(x,y)\le \la^{-1}$ and concludes the proof of \eqref{2.28}.

To prove \eqref{2.28} we can exploit the fact that, unlike the case of $t=0$,
if $t\ne 0$ then $\cos t\sqrt{-\Delta_g}: C^\infty(M)\to C^\infty(M)$ is a conormal
Fourier integral operator with singular support of codimension one.  Based on
this and \eqref{2.17} we deduce that if $(x,t,\xi,\tau; y,\eta)$ is in the
wave front set of 
$$(\cos (t\sqrt{-\Delta_{\tilde g}}))(\tilde x, \alpha_j(\tilde y)),
\quad  j\ne 0, $$
and both $x$ and $y$ are on $\gamma_0$ then both $\xi$ and $\eta$ must be
on the first coordinate axis.  Therefore, since the symbol, $b_\varepsilon(x,\xi)$,
of $b_\varepsilon$
 equals one when $x\in \gamma_0$ and $\xi$ is in a conic neighborhood
of this axis (depending on $\varepsilon$), we conclude that there must be a neighborhood
${\mathcal N}$ of $\gamma_0$ in $M$ so that
$$
\bigl(\, b_\varepsilon \circ \cos (t\sqrt{-\Delta_{\tilde g}})\circ b_\varepsilon^*\,
\bigr)(\tilde x, \alpha_j(\tilde y))
-\bigl(\, \cos t\sqrt{-\Delta_{\tilde g}}\, \bigr)(\tilde x, \alpha_j(\tilde y))
\in C^\infty({\mathcal N}\times {\mathcal N}), \, \, 0\ne |j|\ell(\gamma_0)\le T.$$
%
Because of this, we would have the remaining inequality, \eqref{2.27}, if we could show that
\begin{multline}\label{2.32}
\Bigl| \, \int \beta((s-j\ell(\gamma_0))/5) \Hat \rho(s/T) e^{-is\la}
\bigl(\cos s\sqrt{-\Delta_{\tilde g}}\bigr)(\tilde x, \alpha_j(\tilde y))\, ds \, \Bigr|
\\
\le 
C\la^{\frac12} 
  \, \bigl(\, \max\bigl\{d_g(\tilde x, \alpha_j(\tilde y)), 
e^{\kappa d_g(\tilde x, \alpha_j(\tilde y))}
\bigr\}\, \bigr)^{-\frac12}\,
+C_{T}
\quad x,y\in {\mathcal N}, \, \, 0\ne |j|\ell(\gamma_0)\le T.
\end{multline}

To prove this, we shall use the fact that on $(\Rt, \tilde g)$ we can use the Hadamard parametrix
even for large times.  Recall that the Hadamard parametrix says that if we set
$${\mathcal E}_0(t,x)=(2\pi)^{-2}\int_{\Rt} e^{ix\cdot \xi} \cos(t|\xi|)\, d\xi,$$
and define ${\mathcal E}_\nu$, $\nu=1,2,3,\dots$ recursively by
$2
{\mathcal E}_\nu(t,x)= t\int_0^t{\mathcal E}_{\nu-1}(s,x)ds$,
$\nu=1,2,3,\dots,$
then there are functions $w_\nu\in C^\infty(\Rt\times \Rt)$ so that we have
$$\bigl(\cos(t\sqrt{-\Delta_{\tilde g}})(x,y)
=\sum_{\nu=0}^N w_\nu(x,y) \,  {\mathcal E}_\nu(t, \, d_{\tidle g}(x,y))+R_N(t,x,y),
$$ where for $n=2$, $R_N\in L^\infty_{loc}({\mathbb R}\times \Rt\times \Rt)$ if $N\ge10$.  We are
abusing the notation a bit by putting ${\mathcal E}_\nu(t,r)$ equal to the radial
function ${\mathcal E}_\nu(t,x)$ for some $|x|=r$.   The ${\mathcal E}_\nu$, 
$\nu=1,2,3,\dots$, are Fourier integrals of order $-\nu$; for instance,
$${\mathcal E}_1(t,x)=(2\pi)^{-2}\int_{\Rt}e^{ix\cdot \xi}  \, \frac{t\sin t|\xi|}{2|\xi|}\, d\xi.$$

As a result of this, we would have \eqref{2.32} if we could show that
\begin{multline}\label{2.33}
\Bigl|w_0(\tilde x, \alpha_j(\tilde y))\iint 
\beta((s-j\ell(\gamma_0))/5)\Hat \rho(s/T) \, e^{-i\la s}
e^{i(\tilde x-\alpha_j(\tilde y))\cdot \xi} \cos(s|\xi|)\, d\xi \, ds\, \Bigr|
\\
\le 
C\la^{\frac12} 
  \, \bigl(\, \max\bigl\{d_g(\tilde x, \alpha_j(\tilde y)), 
e^{\kappa d_g(\tilde x, \alpha_j(\tilde y))}
\bigr\}\, \bigr)^{-\frac12},
\quad j=1,2,\dots,
\end{multline}
as well as
\begin{equation}\label{2.34}
\Bigl| \, \int \beta((s-j\ell(\gamma_0))/5) \Hat \rho(s/T) e^{-is\la}
\mathcal{E}_\nu(s, d_g(\tilde x, \alpha_j(\tilde y)) \, ds \, \Bigr|
\le C_\nu, \quad 0\ne j\ell(\gamma_0)\le T, \, \, \nu=1,2,3,\dots.
\end{equation}
Here we are using the fact that $|w_\nu( x, y)|\le C_T$ for $|x|, |y|\le T$.

If we repeat the stationary phase argument that was used to prove \eqref{2.28}, we
see that the left side of \eqref{2.33} is dominated by a fixed constant times
$$\la^{\frac12}\,  w_0(\tilde x, \alpha_j(\tilde y)) \,  (d_{\tilde g}(\tilde x, \alpha_j(\tilde y)))^{-\frac12},
$$
and, consequently, we would have \eqref{2.28} if
\begin{equation}\label{2.35}
w_0(\tilde x, \alpha_j(\tilde y))(d_{\tilde g}(\tilde x, \alpha_j(\tilde y)))^{-\frac12}
\le 
C\bigl(\, \max\bigl\{d_g(\tilde x, \alpha_j(\tilde y)), 
e^{\kappa d_g(\tilde x, \alpha_j(\tilde y))}
\bigr\}\, \bigr)^{-\frac12}
\end{equation}
assuming, as above, that the curvature of $M$ is $\le -\kappa^2$, $\kappa\ge0$.  The last
inequality comes from the fact that in geodesic normal coordinates about $x$, we have
$$w_0(x,y)=\bigl(\, \text{det } g_{ij}(y)\, \bigr)^{-\frac14},$$
(see \cite{Berard}, \cite{Had} or \S2.4 in \cite{SHZ}).  If $y$ has geodesic polar coordinates $(t,\theta)$
about $x$, then $t=d_{\tilde g}(x,y)$, and if ${\mathcal A}(t,\theta)$ is as in \eqref{2.12}, we
conclude that $w_0(x,y)=\sqrt{t/{\mathcal A}(t,\theta)}$, and therefore
\eqref{2.35} follows from G\"unther's comparison estimate \eqref{2.14} if 
$-\kappa^2<0$ and \eqref{2.13} if $\kappa=0$.

The second estimate \eqref{2.34} is elementary and left for the reader, who can
check that the terms are actually $O(\la^{\frac12-\nu})$.  (This is also just a 
special case of Lemma~3.5.3  in \cite{SHZ}.)
This completes the proof of \eqref{2.32}, and, hence, that of Theorem~\ref{theorem1}. \qed
\medskip

\newsection{Concluding remarks}  

It is straightforward to see that the proof of Theorem~\ref{theorem1} shows that one can
strengthen our main estimate \eqref{1.1} in a natural way.  Specifically, if $\gamma_0$
is a periodic geodesic of length  $\ell(\gamma_0)$ and if we define the $\delta$-tube
about $\gamma$ to be
$${\mathcal T}_\delta(\gamma_0)=
\{y\in M: \, \text{dist}_g(y,\gamma_0)<\delta\},$$
with $\delta>0$ fixed, then there is a uniform constant $C_\delta$ so that
whenever $\varepsilon>0$ we have for large $\lambda$
\begin{equation}\label{3.1}
\frac1{\ell(\gamma_0)}\int_{\gamma_0}|e_\la|^2 \, ds
\le \varepsilon \lambda^{\frac12}\|e_\lambda\|^2_{L^2({\mathcal T}_\delta(\gamma_0))}
+C_{\gamma_0,\delta,\varepsilon}\|e_\lambda\|_{L^2(M)}^2.
\end{equation}
Thus, \eqref{1.1} essentially lifts to the cylinder ${\mathbb R}^2/G_{\gamma_0}$, with, as
above, $G_{\gamma_0}$, being the stabilizer group for the lift of $\gamma_0$ to 
the universal cover
$({\mathbb R}^2, \tilde g)$.

To prove this, we as before write $I=B_\varepsilon+b_\varepsilon$, with $b_\varepsilon(x,\xi)$
equal to one near $T^*\gamma_0$ but supported in a small conic neighborhood of this set.
Since the analog of \eqref{2.16} is valid, i.e.,
\begin{equation}\label{3.2}
\int_{\gamma_0} \Bigl| \, \frac1{\pi T}\int_{-T/4}^{T/4} \Hat \chi(t/T) e^{-i\la t}B_\varepsilon
\circ \cos(t\sqrt{-\Delta_g})f\, dt \, \Bigr|^2 \, ds
\le C'_{T,\varepsilon,\gamma_0}\|f\|_{L^2(M)}^2, \end{equation}
it suffices to show that
$$\frac1{\ell(\gamma_0)} \int_{\gamma_0}
\left| \, \frac1T \int_{-T/4}^{T/4} \Hat \chi(t/T) e^{-i\la t} \, b_\varepsilon \circ
\cos \bigl(t\sqrt{-\Delta_g}\bigr) e_\lambda \, dt\, \right|^2 \, ds
$$
is dominated by the right side of \eqref{3.1}.  

If $K_\varepsilon(x,s)$, $x\in M$,
$s\in \gamma_0$ denotes the kernel of this operator then, if $\delta>0$ and $T$ are fixed,
it follows that
\begin{equation}\label{3.3}
|K_\varepsilon(x,s)|\le C_{\gamma_0,T,\delta}, \quad
x\notin {\mathcal T}_\delta(\gamma_0),
\end{equation}
provided that $b_\varepsilon$ is supported in a sufficiently small conic neighborhood of
$T^*\gamma_0$.  This is a simple consequence of the fact that when $b_\varepsilon$ is as
above, by \eqref{2.17},
$\bigl(b_\varepsilon \circ \cos t\sqrt{-\Delta_g}\bigr)(x,s)$ is smooth
when $x\notin {\mathcal T}_\delta(\gamma_0)$,    $s\in \gamma_0$ and $|t|\le T$. Since \eqref{2.20} is valid, we conclude that
there is a uniform constant $C$ so that for large $\lambda$ we have
\begin{multline}\label{3.4}
\frac1{\ell(\gamma_0)} \int_{\gamma_0}
\left| \, \frac1T \int_{-T/4}^{T/4} \Hat \chi(t/T) e^{-i\la t} \, b_\varepsilon \circ
\cos \bigl(t\sqrt{-\Delta_g}\bigr) e_\lambda \, dt\, \right|^2 \, ds
\\
\le CT^{-1}\lambda^{\frac12}\|e_\lambda\|_{L^2({\mathcal T}_\delta(\gamma_0))}^2
+C_{T,\delta,\gamma_0}\|e_\lambda\|_{L^2(M)}^2,
\end{multline}
which along with \eqref{3.2} gives us \eqref{3.1}.  This is because we can dominate
the quantity in \eqref{3.4} by the sum of the corresponding expression where
$e_\lambda$ is replaced by ${\mathbf 1}_{{\mathcal T}_\delta(\gamma_0)}e_\lambda$
and ${\mathbf 1}_{{\mathcal T}^c_\delta(\gamma_0)}e_\lambda$ and use
\eqref{2.20} and our earlier arguments to show that the first of these terms is dominated by the first
term in the right side of \eqref{3.4} if $\lambda$ is large, while the second such term
is dominated by last term in the right side of \eqref{3.4} on account of \eqref{3.3}.

We would also like to point out that
it seems likely that one should be able to take
the parameter $T$ in the proof
of either \eqref{1.1} or \eqref{3.1}
 to be a function of $\lambda$.  This would also require
that the parameter $\varepsilon$ to also be a function of $\lambda$, and thus the argument
would be more involved.  It would not be surprising if, as in B\'erard~\cite{Berard}
or Hassell and Tacey~\cite{HT}, one could
take $T$ to be $\approx \log\la$, in which case the $L^2$-restriction bounds in Theorem~\ref{theorem1} and the $L^4$-estimates in Corollary~\ref{corollary2} could also be improved
to be $O(\la^{\frac14}(\log \la)^{-\delta_1})$ and $O(\la^{\frac18}(\log \la)^{-\delta_2})$,
respectively, for some $\delta_j>0$.  
It is doubtful that these bounds would be optimal, though--indeed if a difficult conjecture
of Rudnick and Sarnak~\cite{sarnak} were valid, both would be $O(\la^\varepsilon)$ for
any $\varepsilon>0$.
One of the main technical issues in carrying out the analysis when $T$ depends on $\la$
would be to determine the analog of \eqref{2.15} in this case.  One would also have to 
take into account more carefully size estimates for the coefficients $w_\nu$, $\nu>0$,
in the Hadamard parametrix, but B\'erard~\cite{Berard} carried out an analysis of these
that would seem to be sufficient if $T\approx \log \la$.   On the other hand, we have argued
here that the $w_0$ coefficient is very well behaved, and so perhaps there could be
further grounds for improvement.  


\end{document}